\documentclass[12pt]{amsart}

\usepackage{amsmath,amssymb,epic,lscape}

\theoremstyle{definition}

\theoremstyle{remark}

\numberwithin{equation}{section}
\providecommand{\bysame}{\leavevmode\hbox to3em{\hrulefill}\thinspace}

\def\DJ{{\hbox{D\kern-.8em\raise.15ex\hbox{--}\kern.35em}}}
\def\DJo{$\;$\kern-.4em
    \hbox{D\kern-.8em\raise.15ex\hbox{--}\kern.35em okovi\'c}}

\def\NSERC{Supported in part by an NSERC Discovery Grant.}

\def\Zn{{\bf Z}_n}

\def\al{{\alpha}}

\def\vf{{\varphi}}
\def\si{{\sigma}}
\def\la{{\lambda}}
\def\bR{{\mbox{\bf R}}}
\def\bZ{{\mbox{\bf Z}}}

\def\pA{{\mathcal A}}
\def\pF{{\mathcal F}}

\renewcommand{\subjclassname}{\textup{2000} Mathematics Subject
Classification }

\begin{document}

\title[Supplementary difference sets with symmetry]
{Supplementary difference sets with symmetry for Hadamard matrices}

\author[D.\v{Z}. \DJ okovi\'{c}]
{Dragomir \v{Z}. \DJ okovi\'{c}}

\address{Department of Pure Mathematics, University of Waterloo,
Waterloo, Ontario, N2L 3G1, Canada}

\email{djokovic@uwaterloo.ca}

\thanks{\NSERC}


\date{}

\begin{abstract}
An overview of the known supplementary difference sets (SDSs) 
$(A_i)$, $1\le i\le4$, with parameters
$(n;k_i;\la)$, $k_i=|A_i|$, where each $A_i$ is either 
symmetric or skew and $\sum k_i=n+\la$ is given.
Five new Williamson matrices over the elementary 
abelian groups of order $5^2$, $3^3$ and $7^2$ are constructed.
New examples of skew Hadamard matrices of order $4n$ for 
$n=47,61,127$ are presented. The last of these is obtained from
a $(127,57,76)$ difference family that we have constructed.
An old non-published example of G-matrices of order 37
is also included.
\end{abstract}

\maketitle
\subjclassname{ 05B20, 05B30 }
\vskip5mm

\section{Introduction}

The Williamson matrices (over any finite abelian group) are too sparse 
\cite{DZ2} to generate Hadamard matrices of all feasible orders. 
The recent extensive computations performed in \cite{HKT} have extended 
the exhaustive searches for circulant Williamson matrices of odd order 
$n$ to  the range $n\le59$. This result 
was made possible not only by using the new and faster computing
devices but also by designing a new efficient algorithm. 
However, the search produced only one new set of Williamson matrices. 
The authors suggest that the researchers should study instead the 
class of Williamson-type matrices. The Williamson matrices over 
arbitrary finite abelian groups belong to this wider class.

In the present paper we consider the method due to Goethals and Seidel 
\cite{GS} of constructing Hadamard matrices by using their well-known 
array
$$
\left[ \begin{array}{rrrr}
		U & XR & YR & ZR \\
		-XR & U & -Z^TR & Y^TR \\
		-YR & Z^TR & U & -X^TR \\
		-ZR & -Y^TR & X^TR & U
\end{array} \right]. $$
One has to find suitable quadruples of $n\times n$ binary 
matrices which can be substituted for $U,X,Y,Z$ in this array 
to give a Hadamard matrix $H$ of order $4n$. (For the symbol $R$ see
the next section.) One way of producing such suitable 
quadruples is via the supplementary difference sets (SDS)
in a finite abelian group $\pA$ of order $n$. 
The parameters of suitable SDSs $A=(A_i)$,  $1\le i\le4$,  must
satisfy an additional condition (see (\ref{uslov-1})). 
If all $A_i$s are symmetric in the sense that $A_i=-A_i$ then
the type I matrices constructed from $A$ are Williamson
matrices in the wider sense. The classical Williamson matrices arise 
when $\pA$ is cyclic. If we only require that the first subset
$A_1$ be skew, in the sense that $\pA$ is the disjoint union
of $A_1$, $-A_1$ and $\{0\}$, then the resulting Hadamard matrix $H$
will be of skew type, i.e., $H-I_{4n}$ is skew-symmetric.
We are mainly interested in the cases where each $A_i$
is either symmetric or skew and we introduce the notion
of (symmetry) types. For instance, the type (ksss)
means that we require $A_1$ to be skew and the other three
$A_i$s to be symmetric. The SDSs for Williamson matrices must
have type (ssss). There are essentially only 
four symmetry types (ssss), (ksss), (kkss) and (kkks),
disregarding the cases where the symmetry is only partial.
The matrices arising from the SDSs having one of these symmetry 
types have been studied for some time by many researchers.
We summarize in Tables 1 and 2 what is known about them for 
small odd values of $n\le 63$.

The new results that we have obtained are presented in the last
section. In particular, we have constructed five new 
multicirculant Williamson matrices, two for each of the orders
25, 27 and one for 49. We also give a new set of $G$-matrices
of order $37$ and new skew Hadamard matrices of order 
$4n$ for $n=47,61,127$. The last one is constructed via 
the new difference family with parameters $(127,57,76)$.
This family also gives a BIBD with the same parameters.

\section{Preliminaries}

Let $\pA$ be a finite abelian group of order $n$. 
Let $A=(A_1,A_2,A_3,A_4)$, $A_i\subseteq \pA$,
be an SDS and let $k_i=|A_i|$ be the cardinality of $A_i$. 
By the definition of SDSs, there exists an integer
$\la\ge0$ such that each nonzero element $a\in\pA$ can be written in
exactly $\la$ ways as the difference $a=x-y$ 
with $\{x,y\}\subseteq A_k$ and $k\in\{1,2,3,4\}$.
We refer to the 6-tuple $(n;k_1,k_2,k_3,k_4;\la)$ as the 
{\em set of parameters} of $A$. We shall be interested only 
in the case when the parameters satisfy the condition
\begin{equation} \label{uslov-1}
 \la=k_1+k_2+k_3+k_4-n. 
\end{equation}
The set of all such SDSs will be denoted by $\pF_\pA$ or just $\pF$.

Let $X=(X_{x,y})$ be an $n\times n$ matrix  
whose rows and columns are indexed by the elements $x,y\in\pA$.
Such $X$ is {\em type I} resp. {\em type II} matrix (relative to 
$\pA$) if $X_{x+z,y+z}=X_{x,y}$ resp. $X_{x+z,y-z}=X_{x,y}$ for all
$x,y,z\in\pA$. 
Let $R$ be the type II matrix defined by $R_{x,y}=\delta_{x+y,0}$,
where $\delta$ is the Kronecker symbol. Then $R^2=I$, the identity 
matrix. The following facts are well known and easy to verify
(see e.g. \cite[Section 1.2]{JS1}).
Any two type I matrices commute. Any type II matrix is symmetric. 
If $X$ and $Y$ are both type I or both type II, then $XY$ is type I.
If $X$ resp. $Y$ is a type I resp. type II matrix, then $XY$ and 
$YX$ are type II, and $X$ and $Y$ are {\em amicable}, i.e., 
$XY^T=YX^T$, where $T$ denotes transposition.
If $X$ and $Y$ are type I  and symmetric, then $XR$ and $YR$ are
amicable and commute.

We say that a matrix is {\em binary} if its entries are $\pm1$.
Let $X\subseteq\pA$ and let $\chi:\pA\to\bR$ be the characteristic 
function of $X$. We denote by $X^c$ the type I binary matrix with 
entries
$$ X^c_{x,y}=1-2\chi(y-x), \quad x,y\in\pA. $$
Thus $X^c_{0,y}=-1$ if and only if $y\in X$.
It is well known that for any $A\in\pF$ the
following matrix equation holds
\begin{equation} \label{uslov-2}
 \sum_{i=1}^4 ({A_i}^c)^T {A_i}^c= 4n I_n.
\end{equation}
Each row-sum of ${A_i}^c$ is equal to $a_i=n-2k_i$. It follows
easily from (\ref{uslov-2}) that 
\begin{equation} \label{uslov-3}
 \sum_{i=1}^4 {a_i}^2= 4n.
\end{equation}

For $X\subseteq\pA$ we say that $X$ is {\em symmetric} resp. {\em skew} if
$-X=X$ resp. $X$, $-X$ and $\{0\}$ form a partition of $\pA$.
If there is a skew $X\subseteq\pA$ then $n$ must
be odd and $|X|=(n-1)/2$. Let $\Sigma=\{s,k,\ast\}$ be the set of
three symbols. We refer to a sequence $(\si_1\si_2\si_3\si_4)$
with $\si_i\in\Sigma$ as a {\em symmetry type} (or simply a {\em type}). 
We say that an SDS $A=(A_i)$ has type $(\si_1\si_2\si_3\si_4)$
if, for each $i$, $A_i$ is symmetric resp. skew when $\si_i=s$ resp. 
$\si_i=k$. No condition is imposed on $A_i$ when $\si_i=\ast$.

When $A$ has type (ssss) then the matrices ${A_i}^c$, $i=1,\ldots,4$, 
are known as the {\em Williamson matrices}. These are four symmetric 
type I binary matrices satisfying the equation (\ref{uslov-2}).

When $A$ has type (ksss) then the matrices
${A_1}^c, {A_2}^cR, {A_3}^cR, {A_4}^cR$ are good matrices. 
When $A$ has type (kkss) or (kkks) then the matrices
$({A_i}^c)$ are G-matrices or best matrices, respectively.
For the general definition of good matrices, G-matrices and best 
matrices see \cite{KS}.

The cyclic case, i.e., when $\pA$ is a cyclic group, has been
investigated most thoroughly. We refer to \cite{HKT} for the
up-to-date information on cyclic Williamson matrices, 
including the complete listing of all non-equivalent such matrices 
of odd order $\le59$. See also the survey papers \cite{SY,HKS}. 
For further information on the other three symmetry types of matrices 
the reader should consult the survey paper \cite{KS} and its references.

For any $A=(A_i)\in\pF$, we can plug the matrices $A_i^c$ into
the Goethals--Seidel array to obtain a Hadamard matrix $H$ of 
order $4n$.  More precisely,
we substitute the symbol $R$ with the $n\times n$ type II
matrix $R$ defined above, and substitute the symbols $U,X,Y,Z$ 
with the four type I matrices $A_1^c,A_2^c,A_3^c,A_4^c$ 
(in that order). If $A$ has type (k$\ast$$\ast$$\ast$), i.e., $A_1$ is skew,
then $H$ will be a skew Hadamard matrix.

Apart from the basic case $\pA=\Zn$ we consider 
here also the case of non-cyclic elementary abelian groups 
in their incarnation as the additive group $(F_q,+)$ 
of a finite field $F_q$ of order $q$. We refer to the
latter type of SDSs as the {\em multicirculant SDSs}.

\section{Known results: Cyclic SDSs}

There are only two known infinite series of cyclic Williamson 
matrices. The first, due to Turyn \cite{RT}, gives Williamson 
matrices of order $(q+1)/2$ where $q$ is a prime power 
$\equiv 1 \pmod{4}$. These matrices are listed on Jennifer 
Seberry's homepage \cite{JS2} for orders $\le63$.
The second, due to Whiteman \cite{AW},
gives Williamson matrices of order $p(p+1)/2$ where $p$ is a 
prime $\equiv 1 \pmod{4}$.
There is also an infinite series of cyclic $G$-matrices constructed 
by Spence \cite{ES}. Their orders are $(q+1)/2$ where $q$ is a
prime power $\equiv 5 \pmod{8}$. 
We are not aware of the existence of any infinite series of 
good or best matrices.

In Table 1 we summarize what is known about the existence of
cyclic SDSs $A\in\pF$ with specified symmetry (ssss), (ksss), 
(kkss) or (kkks) for small odd values of $n$ ($\le 63$). 
For the entry of Table 1 (and those of Table 2)
marked with the symbol $\dag$ see Section \ref{NoviRez}.

In the first three columns we list the feasible parameters 
$n$, $(k_i)$, $\la$ with 
$k_1\ge k_2\ge k_3\ge k_4$ and $2k_1<n$. Note that these 
conditions are not restrictive since we can permute the 
$A_i$s and replace any $A_i$ with its complement.
As the row-sums $a_i$ of the matrices ${A_i}^c$ are often used, 
we list them in the fourth column. 
By our choice of the $k_i$ we have $a_i>0$ for all $i$. 
For each of the above four symmetry types we give in
the last four columns the number
of known non-equivalent SDSs. If this number is 
written in bold type then an exhaustive search for these families 
has been carried out and a reference is provided. 
The sign $\times$ means that the parameter set is not compatible 
with the symmetry type of the column, and the blank entry
means that the existence question remains unresolved.
(The second example of G-matrices for $n=41$ 
given in \cite{GK} is not valid.)

In the case of $G$-matrices constructed by Spence one has
$k_1=k_2=(n-1)/2$. This determines uniquely $k_3$ and $k_4$
in the cases $n=51,55$ but not in the case $n=63$. In the last
case we had to construct explicitly the SDS by using linear
recurrent sequences as explained in \cite{ES} and its
references. Since this was quite involved computation, we
sketch here some details.

We start with the finite field $F_q=\bZ_5/(x^3-2x+2)$ of order 
$q=5^3=125$ and denote by $a$ the image of the variable $x$.
The polynomial $x^3-2x+2$ is primitive over $\bZ_5$, 
i.e., $a$ generates the multiplicative group $F^*_q$. 
We consider next the linear reccurence relation 
$ax_{i+1}+x_i+x_{i-1}=0$, $i=1,2,\ldots$, with 
initial values $x_0=x_1=1$. One can verify that the 
infinite sequence $(x_0,x_1,x_2,\ldots)$ 
generated by the above relation has minimal 
period $q^2-1=15624$, i.e., it is an $m$-sequence in the
terminology of \cite{ES}. The set of indexes
$X=\{i:0\le i<q^2-1,\ x_i=1\}$ is a cyclic relative difference
set with parameters $(126,124,125,1)$ using the definition
in \cite{IK,ES}. By reducing these indexes modulo $4(q+1)=504$, 
we obtain the cyclic relative difference set $Y$ with 
parameters $(63,8,125,31)$. By replacing $Y$ with the translate
$Y+113\subseteq\bZ_{504}$, we obtain a $Y$ which is
fixed under multiplication by $q$:
\begin{eqnarray*}
&& Y=\{ 8,9,11,12,16,17,19,21,24,26,38,39,40,41,42,44,45,53, \\
&& 54,55,59,60,62,73,80,81,83,85,91,92,95,96,98,103,104,105, \\
&& 106,109,117,119,120,122,128,130,136,146,154,176,177,183,190, \\
&& 195,198,200,204,205,210,214,220,225,226,237,249,252,253,257, \\
&& 259,265,266,270,275,277,283,284,287,295,300,304,310,313,317, \\
&& 319,322,323,328,339,342,353,359,365,367,368,373,376,377,381, \\
&& 384,393,400,405,407,408,411,412,414,415,424,425,427,434,444, \\
&& 446,453,455,460,464,467,471,475,480,486,488,490,492,496 \}.
\end{eqnarray*}

For $1\le i\le4$ let $Y_i=\{j\in Y:j\equiv i-1 \pmod{8} \}$ and let
$A_i=Y_i \pmod{63}$. The blocks $A_1$ and $A_3$ are symmetric while
$A_2$ and $A_4$ are skew. Thus they are uniquely determined by the
intersections $A_i^*=A_i\cap\{0,1,\ldots,31\}$. Explicitly, we have
\begin{eqnarray*}
A_1^* &=& \{2,4,6,7,8,10,11,13,14,16,17,20,21,22,23,24,26,28,30\}, \\
A_2^* &=& \{ 5,6,7,8,9,10,11,13,15,17,18,19,20,23,24,26,28,31\}, \\
A_3^* &=& \{ 0,3,4,7,12,14,15,19,20,21,26,28,29,31\}, \\
A_4^* &=& \{ 1,2,3,6,7,8,11,12,16,19,20,21,22,23,24,25,26,27,28,31\}.
\end{eqnarray*}
Finally we replace $A_1$ with its complement. After permuting the blocks,
the new SDS has parameters $(63;31,31,27,25;51)$ and type (kkss). \\

\newpage

\begin{center}
\begin{tabular}{|r|c|c|c|c|c|c|c|}
\multicolumn{8}{c}{\bf Table 1: Cyclic SDSs with symmetry} \\
\multicolumn{8}{c}{} \\ \hline 
\multicolumn{1}{|r|}{$n$} & \multicolumn{1}{c|}{ $(k_i)$ } & 
\multicolumn{1}{c|}{ $\la$ } & \multicolumn{1}{c|}{ $(a_i)$ } & 
\multicolumn{1}{c|}{(ssss)} & \multicolumn{1}{c|}{(ksss)} &
\multicolumn{1}{c|}{(kkss)} & \multicolumn{1}{c|}{(kkks)} \\ \hline
3 & 1,1,1,0 & 0 & 1,1,1,3 & {\bf 1},\cite{BH} & {\bf 1},\cite{DH} & 
{\bf 1},\cite{GK} & {\bf 1},\cite{GKJ} \\ \hline
5 & 2,2,1,1 & 1 & 1,1,3,3 & {\bf 1},\cite{BH} & {\bf 1},\cite{DH} & 
{\bf 1},\cite{GK} & $\times$ \\ \hline
7 & 3,3,3,1 & 3 & 1,1,1,5 & {\bf 1},\cite{BH} & {\bf 1},\cite{DH} & 
{\bf 1},\cite{GK} & {\bf 1},\cite{GKJ} \\
  & 3,2,2,2 & 2 & 1,3,3,3 & {\bf 1},\cite{BH} & {\bf 2},\cite{DH} & 
$\times$ & $\times$ \\ \hline
9 & 4,4,3,2 & 4 & 1,1,3,5 & {\bf 2},\cite{BH} & {\bf 1},\cite{DH} & 
{\bf 1},\cite{GK} & $\times$ \\
  & 3,3,3,3 & 3 & 3,3,3,3 & {\bf 1},\cite{BH} & $\times$ & $\times$ & 
$\times$ \\ \hline
11 & 5,4,4,3 & 5 & 1,3,3,5 & {\bf 1},\cite{BH} & {\bf 3},\cite{Sz} & 
$\times$ & $\times$ \\ \hline
13 & 6,6,6,3 & 8 & 1,1,1,7 & {\bf 1},\cite{BH} & {\bf 2},\cite{DH} & 
{\bf 0},\cite{GK} & {\bf 2},\cite{GKJ} \\
   & 6,6,4,4 & 7 & 1,1,5,5 & {\bf 1},\cite{BH} & {\bf 4},\cite{DH} & 
{\bf 8},\cite{GK} & $\times$ \\
   & 5,5,5,4 & 6 & 3,3,3,5 & {\bf 2},\cite{BH} & $\times$ & $\times$ & 
$\times$ \\ \hline
15 & 7,7,6,4 & 9 & 1,1,3,7 & {\bf 3},\cite{BH} & {\bf 7},\cite{DH} & 
{\bf 32},\cite{GK} & $\times$ \\
   & 7,6,5,5 & 8 & 1,3,5,5 & {\bf 1},\cite{BH} & {\bf 4},\cite{DH} & 
$\times$ & $\times$ \\ \hline
17 & 8,7,7,5 & 10 & 1,3,3,7 & {\bf 3},\cite{BH} & {\bf 2},\cite{DH} & 
$\times$ & $\times$ \\
   & 7,7,6,6 & 9 & 3,3,5,5 & {\bf 1},\cite{BH} & $\times$ & $\times$ & 
$\times$ \\ \hline
19 & 9,9,7,6 & 12 & 1,1,5,7 & {\bf 3},\cite{BH} & {\bf 5},\cite{DH} & 
{\bf 9},\cite{GK} & $\times$ \\
   & 8,8,8,6 & 11 & 3,3,3,7 & {\bf 3},\cite{BH} & $\times$ & $\times$ & 
$\times$ \\
   & 9,7,7,7 & 11 & 1,5,5,5 & {\bf 0},\cite{BH} & {\bf 3},\cite{DH} & 
$\times$ & $\times$ \\ \hline
21 & 10,10,10,6 & 15 & 1,1,1,9 & {\bf 1},\cite{BH} & {\bf 4},\cite{DH} & 
{\bf 23},\cite{GK} & {\bf 21},\cite{GKJ} \\
   & 10,9,8,7 & 13 & 1,3,5,7 & {\bf 3},\cite{BH} & {\bf 6},\cite{DH} & 
$\times$ & $\times$ \\
   & 9,8,8,8 & 12 & 3,5,5,5 & {\bf 3},\cite{BH} & $\times$ & $\times$ & 
$\times$ \\ \hline
23 & 11,11,10,7 & 16 & 1,1,3,9 & {\bf 0},\cite{BH} & {\bf 6},\cite{Sz} & 
{\bf 16},\cite{GK} & $\times$ \\
   & 10,10,9,8 & 14 & 3,3,5,7 & {\bf 1},\cite{BH} & $\times$ & $\times$ & 
$\times$ \\ \hline
25 & 12,11,11,8 & 17 & 1,3,3,9 & {\bf 1},\cite{DZ4} & {\bf 3},\cite{Sz} & 
$\times$ & $\times$ \\
   & 12,12,9,9 & 17 & 1,1,7,7 & {\bf 3},\cite{DZ4} & {\bf 0},\cite{Sz} & 
{\bf 13},\cite{GK} & $\times$ \\
   & 12,10,10,9 & 16 & 1,5,5,7 & {\bf 3},\cite{DZ4} & {\bf 6},\cite{Sz} & 
$\times$ & $\times$ \\
   & 10,10,10,10 & 15 & 5,5,5,5 & {\bf 3},\cite{DZ4} & $\times$ & $\times$ & 
$\times$ \\ \hline
27 & 13,13,11,9 & 19 & 1,1,5,9 & {\bf 2},\cite{KSa} & {\bf 6},\cite{Sz} & 
{\bf 20},\cite{GK} & $\times$ \\
   & 12,12,12,9 & 18 & 3,3,3,9 & {\bf 0},\cite{KSa} & $\times$ & $\times$ & 
$\times$ \\
   & 13,12,10,10 & 18 & 1,3,7,7 & {\bf 3},\cite{KSa} & {\bf 6},\cite{Sz} & 
$\times$ & $\times$ \\
   & 12,11,11,10 & 17 & 3,5,5,7 & {\bf 1},\cite{KSa} & $\times$ & $\times$ & 
$\times$ \\ \hline
29 & 14,13,12,10 & 20 & 1,3,5,9 & {\bf 1},\cite{DZ1} & {\bf 5},\cite{Sz} & 
$\times$ & $\times$ \\
   & 13,13,11,11 & 19 & 3,3,7,7 & {\bf 0},\cite{DZ1} & $\times$ & $\times$ & 
$\times$ \\ \hline
31 & 15,15,15,10 & 24 & 1,1,1,11 & {\bf 0},\cite{DZ1} & {\bf 2},\cite{Sz} & 
{\bf 8},\cite{GK} & {\bf 8},\cite{GKJ} \\
   & 14,14,13,11 & 21 & 3,3,5,9 & {\bf 0},\cite{DZ1} & $\times$ & 
$\times$ & $\times$ \\
   & 15,13,12,12 & 21 & 1,5,7,7 & {\bf 1},\cite{DZ1} & {\bf 1},\cite{Sz} & 
$\times$ & $\times$ \\
   & 13,13,13,12 & 20 & 5,5,5,7 & {\bf 1},\cite{DZ1} & $\times$ & $\times$ & 
$\times$ \\ \hline
\end{tabular} \\
\end{center}

\begin{center}
\begin{tabular}{|r|c|c|c|c|c|c|c|}
\multicolumn{8}{c}{\bf Table 1 (continued)} \\ \hline 
\multicolumn{1}{|r|}{$n$} & \multicolumn{1}{c|}{ $(k_i)$ } & 
\multicolumn{1}{c|}{ $\la$ } & \multicolumn{1}{c|}{ $(a_i)$ } & 
\multicolumn{1}{c|}{(ssss)} & \multicolumn{1}{c|}{(ksss)} &
\multicolumn{1}{c|}{(kkss)} & \multicolumn{1}{c|}{(kkks)} \\ \hline
33 & 16,16,15,11 & 25 & 1,1,3,11 & {\bf 1},\cite{DZ2} & {\bf 6},\cite{GKS} & 
{\bf 9},\cite{GK} & $\times$ \\
   & 16,16,13,12 & 24 & 1,1,7,9 & {\bf 1},\cite{DZ2} & {\bf 4},\cite{GKS} & 
{\bf 22},\cite{GK} & $\times$ \\
   & 16,14,14,12 & 23 & 1,5,5,9 & {\bf 2},\cite{DZ2} & {\bf 5},\cite{GKS} & 
$\times$ & $\times$ \\
   & 15,14,13,13 & 22 & 3,5,7,7 & {\bf 1},\cite{DZ2} & $\times$ & 
$\times$ & $\times$ \\ \hline
35 & 17,16,16,12 & 26 & 1,3,3,11 & {\bf 0},\cite{DZ2} & {\bf 4},\cite{GKS} & 
$\times$ & $\times$ \\
   & 17,16,14,13 & 25 & 1,3,7,9 & {\bf 0},\cite{DZ2} & {\bf 2},\cite{GKS} & 
$\times$ & $\times$ \\
   & 16,15,15,13 & 24 & 3,5,5,9 & {\bf 0},\cite{DZ2} & $\times$ & $\times$ & 
$\times$ \\ \hline
37 & 18,18,16,13 & 28 & 1,1,5,11 & {\bf 0},\cite{DZ4} & {\bf 1},\cite{GKS} & 
5,\cite{GK},$\dag$ & $\times$ \\
   & 17,17,17,13 & 27 & 3,3,3,11 & {\bf 1},\cite{DZ4} & $\times$ & $\times$ & 
$\times$ \\
   & 17,17,15,14 & 26 & 3,3,7,9 & {\bf 1},\cite{DZ4} & $\times$ & $\times$ & 
$\times$ \\
   & 18,15,15,15 & 26 & 1,7,7,7 & {\bf 0},\cite{DZ4} & {\bf 1},\cite{GKS} & 
$\times$ & $\times$ \\
   & 16,16,15,15 & 25 & 5,5,7,7 & {\bf 2},\cite{DZ4} & $\times$ & $\times$ & 
$\times$ \\ \hline
39 & 19,18,17,14 & 29 & 1,3,5,11 & {\bf 0},\cite{DZ2} & {\bf 3},\cite{GKS} & 
$\times$ & $\times$ \\
   & 19,17,16,15 & 28 & 1,5,7,9 & {\bf 0},\cite{DZ2} & {\bf 2},\cite{GKS} & 
$\times$ & $\times$ \\
   & 17,17,17,15 & 27 & 5,5,5,9 & {\bf 1},\cite{DZ2} & $\times$ & $\times$ & 
$\times$ \\
   & 18,16,16,16 & 27 & 3,7,7,7 & {\bf 0},\cite{DZ2} & $\times$ & $\times$ & 
$\times$ \\ \hline
41 & 19,19,18,15 & 30 & 3,3,5,11 & {\bf 0},\cite{HKT} & $\times$ & $\times$ & 
$\times$ \\
   & 20,20,16,16 & 31 & 1,1,9,9 & {\bf 1},\cite{HKT} &  & 1,\cite{GK} & 
$\times$ \\
   & 19,18,17,16 & 29 & 3,5,7,9 & {\bf 0},\cite{HKT} & $\times$ & $\times$ & 
$\times$ \\ \hline
43 & 21,21,21,15 & 35 & 1,1,1,13 & {\bf 0},\cite{HKT} &  &  &  \\
   & 21,21,18,16 & 33 & 1,1,7,11 & {\bf 0},\cite{HKT} &  &  & $\times$ \\
   & 21,19,19,16 & 32 & 1,5,5,11 & {\bf 1},\cite{HKT} &  & $\times$ & 
$\times$ \\
   & 21,20,17,17 & 32 & 1,3,9,9 & {\bf 0},\cite{HKT} &  & $\times$ & $\times$ \\
   & 19,18,18,18 & 30 & 5,7,7,7 & {\bf 1},\cite{HKT} & $\times$ & $\times$ & 
$\times$ \\ \hline
45 & 22,22,21,16 & 36 & 1,1,3,13 & {\bf 0},\cite{RV} &  &  & $\times$ \\
   & 22,21,19,17 & 34 & 1,3,7,11 & {\bf 0},\cite{RV} & $\times$ & $\times$ & 
$\times$ \\
   & 21,20,20,17 & 33 & 3,5,5,11 & {\bf 0},\cite{RV} & $\times$ & $\times$ & 
$\times$ \\
   & 21,21,18,18 & 33 & 3,3,9,9 & {\bf 0},\cite{RV} & $\times$ & $\times$ & 
$\times$ \\
   & 22,19,19,18 & 33 & 1,7,7,9 & {\bf 0},\cite{RV} &  & $\times$ & $\times$ \\
   & 20,20,19,18 & 32 & 5,5,7,9 & {\bf 1},\cite{RV} & $\times$ & $\times$ & 
$\times$ \\ \hline
47 & 23,22,22,17 & 37 & 1,3,3,13 & {\bf 0},\cite{HKT} &  & $\times$ & 
$\times$ \\
   & 22,22,20,18 & 35 & 3,3,7,11 & {\bf 0},\cite{HKT} & $\times$ & 
$\times$ & $\times$ \\
   & 23,21,19,19 & 35 & 1,5,9,9 & {\bf 0},\cite{HKT} &  & $\times$ & $\times$ \\
   & 22,20,20,19 & 34 & 3,7,7,9 & {\bf 0},\cite{HKT} & $\times$ & $\times$ & 
$\times$ \\ \hline
49 & 24,24,22,18 & 39 & 1,1,5,13 & {\bf 0},\cite{HKT} &  &  & $\times$ \\
   & 23,23,23,18 & 38 & 3,3,3,13 & {\bf 0},\cite{HKT} & $\times$ & $\times$ & 
$\times$ \\
   & 24,22,21,19 & 37 & 1,5,7,11 & {\bf 0},\cite{HKT} &  & $\times$ & 
$\times$ \\
   & 22,22,22,19 & 36 & 5,5,5,11 & {\bf 0},\cite{HKT} & $\times$ & $\times$ & 
$\times$ \\
   & 23,22,20,20 & 36 & 3,5,9,9 & {\bf 1},\cite{HKT} & $\times$ & $\times$ & 
$\times$ \\
   & 21,21,21,21 & 35 & 7,7,7,7 & {\bf 0},\cite{HKT} & $\times$ & $\times$ & 
$\times$ \\ \hline
\end{tabular} \\
\end{center}

\begin{center}
\begin{tabular}{|r|c|c|c|c|c|c|c|}
\multicolumn{8}{c}{\bf Table 1 (continued)} \\ \hline 
\multicolumn{1}{|r|}{$n$} & \multicolumn{1}{c|}{ $(k_i)$ } & 
\multicolumn{1}{c|}{ $\la$ } & \multicolumn{1}{c|}{ $(a_i)$ } & 
\multicolumn{1}{c|}{(ssss)} & \multicolumn{1}{c|}{(ksss)} &
\multicolumn{1}{c|}{(kkss)} & \multicolumn{1}{c|}{(kkks)} \\ \hline
51 & 25,24,23,19 & 40 & 1,3,5,13 & {\bf 0},\cite{RV} &  & $\times$ & 
$\times$ \\
   & 25,25,21,20 & 40 & 1,1,9,11 & {\bf 1},\cite{RV} &  & 1,\cite{ES} & $\times$ \\
   & 24,23,22,20 & 38 & 3,5,7,11 & {\bf 1},\cite{RV} & $\times$ & 
$\times$ & $\times$ \\
   & 23,22,22,21 & 37 & 5,7,7,9 & {\bf 0},\cite{RV} & $\times$ & 
$\times$ & $\times$ \\ \hline
53 & 25,25,24,20 & 41 & 3,3,5,13 & {\bf 0},\cite{HKT} & $\times$ & 
$\times$ & $\times$ \\
   & 26,25,22,21 & 41 & 1,3,9,11 & {\bf 0},\cite{HKT} &  & $\times$ & 
$\times$ \\
   & 26,23,22,22 & 40 & 1,7,9,9 & {\bf 0},\cite{HKT} &  & $\times$ & 
$\times$ \\
   & 24,24,22,22 & 39 & 5,5,9,9 & {\bf 0},\cite{HKT} & $\times$ & 
$\times$ & $\times$ \\ \hline
55 & 27,27,24,21 & 44 & 1,1,7,13 & {\bf 0},\cite{HKT} &  & 1,\cite{ES} & $\times$ \\
   & 27,25,25,21 & 43 & 1,5,5,13 & {\bf 0},\cite{HKT} &  & $\times$ & 
$\times$ \\
   & 26,26,23,22 & 42 & 3,3,9,11 & {\bf 1},\cite{HKT} & $\times$ & 
$\times$ & $\times$ \\
   & 27,24,24,23 & 43 & 1,7,7,11 & {\bf 0},\cite{HKT} &  & $\times$ & 
$\times$ \\
   & 25,25,24,22 & 41 & 5,5,7,11 & {\bf 0},\cite{HKT} & $\times$ & 
$\times$ & $\times$ \\
   & 26,24,23,23 & 41 & 3,7,9,9 & {\bf 0},\cite{HKT} & $\times$ & 
$\times$ & $\times$ \\ \hline
57 & 28,28,28,21 & 48 & 1,1,1,15 & {\bf 0},\cite{HKT} &  &  &  \\
   & 28,27,25,22 & 45 & 1,3,7,13 & {\bf 0},\cite{HKT} &  & $\times$ & 
$\times$ \\
   & 27,26,26,22 & 44 & 3,5,5,13 & {\bf 0},\cite{HKT} & $\times$ & 
$\times$ & $\times$ \\
   & 28,26,24,23 & 44 & 1,5,9,11 & {\bf 0},\cite{HKT} &  & $\times$ & 
$\times$ \\
   & 27,25,25,23 & 44 & 3,7,7,11 & {\bf 0},\cite{HKT} & $\times$ & 
$\times$ & $\times$ \\
   & 25,25,25,24 & 42 & 7,7,7,9 & {\bf 1},\cite{HKT} & $\times$ & 
$\times$ & $\times$ \\ \hline
59 & 29,29,28,22 & 49 & 1,1,3,15 & {\bf 0},\cite{HKT} &  &  & $\times$ \\
   & 28,28,26,23 & 46 & 3,3,7,13 & {\bf 0},\cite{HKT} & $\times$ & 
$\times$ & $\times$ \\
   & 28,27,25,24 & 45 & 3,5,9,11 & {\bf 0},\cite{HKT} & $\times$ & 
$\times$ & $\times$ \\
   & 27,26,25,25 & 44 & 5,7,9,9 & {\bf 0},\cite{HKT} & $\times$ & 
$\times$ & $\times$ \\ \hline
61 & 30,29,29,23 & 50 & 1,3,3,15 &  &  & $\times$ & $\times$ \\
   & 30,28,27,24 & 48 & 1,5,7,13 &  &  & $\times$ & $\times$ \\
   & 28,28,28,24 & 47 & 5,5,5,13 &  & $\times$ & $\times$ & $\times$ \\
   & 30,30,25,25 & 49 & 1,1,11,11 & 1,\cite{RT} &  &  & $\times$ \\
   & 28,27,27,25 & 46 & 5,7,7,11 &  & $\times$ & $\times$ & $\times$ \\
   & 30,26,26,26 & 47 & 1,9,9,9 &  &  & $\times$ & $\times$ \\ \hline
63 & 31,31,29,24 & 52 & 1,1,5,15 &  &  &  & $\times$ \\
   & 30,30,30,24 & 51 & 3,3,3,15 &  & $\times$ & $\times$ & $\times$ \\
   & 31,31,27,25 & 51 & 1,1,9,13 &  &  & 1,\cite{ES} & $\times$ \\
   & 30,29,28,25 & 49 & 3,5,7,13 &  & $\times$ & $\times$ & $\times$ \\
   & 31,30,26,26 & 50 & 1,3,11,11 & 1,\cite{RT} &  & $\times$ & $\times$ \\
   & 31,28,27,26 & 49 & 1,7,9,11 &  &  & $\times$ & $\times$ \\
   & 29,29,27,26 & 48 & 5,5,9,11 &  & $\times$ & $\times$ & $\times$ \\
   & 30,27,27,27 & 48 & 3,9,9,9 &  & $\times$ & $\times$ & 
$\times$ \\ \hline
\end{tabular} \\
\end{center}

\section{Known results: Multicirculant SDSs}

There is an infinite series of multicirculant Williamson matrices
due to Xia and Liu \cite{XL1}. It gives matrices of order 
$q^2$ where $q$ is a prime power $\equiv 1 \pmod{4}$. 
For such $q$ they construct SDSs having symmetry type (ssss)
and parameters
$$ \left( q^2; \binom{q}{2}, \binom{q}{2}, \binom{q}{2}, 
\binom{q}{2}; q(q-2) \right). $$

There are only four proper odd prime powers: $3^2,5^2,3^3,7^2$ 
in the range that we consider. If $n$ is one of these powers then
$4n-3$ is not a square. Thus, the symmetry type (kkks) cannot occur.
Table 2 shows what is presently known about the existence of 
multicirculant SDSs for these four powers. It includes the four
previously known isolated examples. The ``No'' entry means that we 
have carried out an exhaustive search and did not find any SDSs 
of that type. \\

\begin{center}
\begin{tabular}{|r|c|c|c|c|c|c|}
\multicolumn{7}{c}{\bf Table 2: Multicirculant SDSs with symmetry} \\
\multicolumn{7}{c}{} \\ \hline 
\multicolumn{1}{|r|}{$n$} & \multicolumn{1}{c|}{ $(k_i)$ } & 
\multicolumn{1}{c|}{ $\la$ } & \multicolumn{1}{c|}{ $(a_i)$ } & 
\multicolumn{1}{c|}{(ssss)} & \multicolumn{1}{c|}{(ksss)} &
\multicolumn{1}{c|}{(kkss)} \\ \hline
$3^2$ & 4,4,3,2 & 4 & 1,1,3,5 & No & No & No \\
  & 3,3,3,3 & 3 & 3,3,3,3 & Yes \cite{SY} & $\times$ & $\times$ \\ \hline
$5^2$ & 12,11,11,8 & 17 & 1,3,3,9 & Yes $\dag$ & No & $\times$ \\
   & 12,12,9,9 & 17 & 1,1,7,7 & Yes $\dag$ & No & Yes \cite{DZ3} \\
   & 12,10,10,9 & 16 & 1,5,5,7 & No & No & $\times$ \\
   & 10,10,10,10 & 15 & 5,5,5,5 & Yes \cite{SY,XL1} & $\times$ & 
$\times$ \\ \hline
$3^3$ & 13,13,11,9 & 19 & 1,1,5,9 & No & No & No \\
   & 12,12,12,9 & 18 & 3,3,3,9 & Yes $\dag$ & $\times$ & $\times$ \\
   & 13,12,10,10 & 18 & 1,3,7,7 & No & No & $\times$ \\
   & 12,11,11,10 & 17 & 3,5,5,7 & No & $\times$ & $\times$ \\ \hline
$7^2$ & 24,24,22,18 & 39 & 1,1,5,13 &  &  & \\
   & 23,23,23,18 & 38 & 3,3,3,13 &  & $\times$ & $\times$ \\
   & 24,22,21,19 & 37 & 1,5,7,11 &  &  & $\times$ \\
   & 22,22,22,19 & 36 & 5,5,5,11 &  & $\times$ & $\times$ \\
   & 23,22,20,20 & 36 & 3,5,9,9 &  & $\times$ & $\times$ \\
   & 21,21,21,21 & 35 & 7,7,7,7 & Yes \cite{QX},$\dag$ & $\times$ & 
$\times$ \\ \hline
\end{tabular} \\
\end{center}

\section{ New results} \label{NoviRez}

The new results of positive nature are presented in increasing order $n$ 
of the additive abelian group $\pA$ employed. 

\subsection{Multicirculant Williamson matrices of order $25$}
Let $F_{25}=\bZ_5[x]/(x^2+2)$ be the finite field of order 25, 
and let us identify $x$ with its image in $F_{25}$. Let
\begin{eqnarray*}
A'_1 &=& \{1,2,x,1+x,2+x,2+2x\}, \\
A'_2 &=& \{x,1+x,1-2x,2\pm x\}, \\
A'_3 &=& \{1,x,2x,2+x,1+2x\}, \\
A'_4 &=& \{1,1\pm x,2-2x\}, \\
B'_1 &=& \{1,2,1+x,1-2x,2+x,2+2x\}, \\
B'_2 &=& \{1,2,x,2x,1+x,2+2x\}, \\
B'_3 &=& \{1,2,2-x,2-2x\}, \\
B'_4 &=& \{1\pm x,1+2x,2+2x\}.
\end{eqnarray*}
The eight subsets
\begin{eqnarray*}
A_1 &=& A'_1\cup(-A'_1),\quad A_2=A'_2\cup \{0\} \cup(-A'_2), \\
A_3 &=& A'_3\cup \{0\} \cup(-A'_3),\quad A_4=A'_4\cup(-A'_4), \\
B_1 &=& B'_1\cup(-B'_1),\quad B_2=B'_2\cup(-B'_2), \\
B_3 &=& B'_3\cup \{0\} \cup(-B'_3),\quad B_4=B'_4\cup \{0\} \cup(-B'_4)
\end{eqnarray*}
are obviously symmetric. One can easily verify that 
$(A_i)$ and $(B_i)$ are SDSs in $\pA=(F_{25},+)$.
Their parameters are $(5^2;12,11,11,8;17)$ and
$(5^2;12,12,9,9;17)$, respectively. 

As far as we know, the existence of elementary abelian SDSs of
type (ssss) and with the above parameters 
was not known previously. For the parameters 
$(5^2;10,10,10,10;15)$ such SDS was constructed by A. Whiteman, 
see \cite{SY}. It turns out that his SDS is equivalent to one 
in the infinite series of M. Xia and G. Liu \cite{XL1}.

\subsection{Multicirculant Williamson matrices of order $27$}
\label{Slucaj27}
Let $F_{27}=\bZ_3[x]/(x^3-x+1)$, a finite field of order 27, 
and let us identify $x$ with its image.
As far as we know, the existence of an elementary abelian 
SDS with parameters $(3^3;12,12,12,9;18)$ and symmetry 
type (ssss) is not known. In the cyclic case it is known \cite{KS} 
that such SDS does not exist. 
We have constructed the following two non-equivalent examples of 
multicirculant SDSs with the above parameters and type. 

Let us begin with the seven subsets
\begin{eqnarray*}
A'_1 &=& \{ 1,x^2,1+x^2,x\pm x^2,1-x-x^2 \}, \\
A'_2 &=& \{ 1,1+x^2,x\pm x^2,1\pm x-x^2 \}, \\
A'_3 &=& \{ 1,x,x^2,1-x^2,x-x^2,1+x-x^2 \}, \\
B'_1 &=& \{ 1,x,1+x,x+x^2,1\pm x+x^2 \}, \\
B'_2 &=& \{ x,x^2,1+x,1-x^2,x+x^2,1+x-x^2 \}, \\
B'_3 &=& \{ x,x^2,x-x^2,1+x+x^2,x^2-x\pm1 \}, \\
A'_4 &=& B'_4=\{ 1,x,1-x^2,x-x^2 \}
\end{eqnarray*}
of $\pA=(F_{27},+)$. Each of the subsets
$$ A_i=A'_i\cup(-A'_i),\quad B_i=B'_i\cup(-B'_i),
\quad i=1,2,3; $$
and also $A_4=B_4=A'_4\cup \{0\} \cup(-A'_4)$
is symmetric. Moreover one can verify that each of the
quadruples $(A_i)$ and $(B_i)$ is an SDS with the
above parameters. Hence the corresponding multicirculant 
matrices, i.e., type 1 matrices, are Williamson matrices.

Let us prove that these two SDSs are not equivalent.
Assume that $\vf(B_i)=A_1+a$ for some automorphism $\vf$
of $\pA$, some $i\in\{1,2,3\}$, and some nonzero element 
$a\in\pA$. Since $B_i=-B_i$ and $A_1=-A_1$, we have
$$ A_1+a=\vf(-B_i)=-\vf(B_i)=-(A_1+a)=A_1-a, $$
and so $A_1=A_1-2a=A_1+a$. This means that $A_1$ is the union of
four cosets of the subgroup $\{0,a,-a\}$. Consequently,
$a$ must occur exactly 12 times in the list of differences
$x-y$ with $x,y\in A_1$. Since $a\ne0$, a simple computation
shows that this is not true. We now conclude that if
our two SDSs are equivalent then there exists an automorphism
$\vf$ of $\pA$ and a permutation $\sigma$ of $\{1,2,3\}$
such that $\vf(A_i)=B_{\sigma(i)}$ for $i=1,2,3$.
Since $|A_1\cap A_2\cap A_3|=4$ and
$|B_1\cap B_2\cap B_3|=2$, this is impossible.
Hence the two SDSs are not equivalent.

\subsection{New $G$-matrices of order $37$}
For $n=37$ four non-equivalent SDSs of type (kkss) were found
in \cite{GK}. (Their search in this case was not exhaustive.)
We have constructed one such SDS in 1995 but were not able to 
include it in our paper \cite{DZ3} and so it remained unpublished.
As it is not equivalent to the four SDSs just mentioned, 
we list it here:
\begin{eqnarray*}
&& (37;18,18,16,13;28) \\
&& \{2,3,5,6,9,10,11,13,15,18,20,21,23,25,29,30,33,36\}, \\
&& \{1,2,4,6,9,10,11,12,17,18,21,22,23,24,29,30,32,34\}, \\
&& \{1,2,4,5,6,10,17,18,19,20,27,31,32,33,35,36\}, \\
&& \{0,3,11,13,15,16,17,20,21,22,24,26,34\}.
\end{eqnarray*}

\subsection{A new skew Hadamard matrix of order $4\cdot47$}
We have constructed recently \cite{DZ5} SDSs with parameters 
$(47;30,22,22;39)$ and $(47;21,19,19;24)$ (two of each kind). 
By combining them with the skew cyclic $(47;23;11)$ difference set, 
we obtained SDSs with parameters $(47;23,30,22,22;50)$ and 
$(47;23,21,19,19;35)$. By replacing in the former the second
set with its complement, the parameters become
$(47;23,22,22,17;37)$. All of these SDSs have
symmetry type (k$\ast$$\ast$$\ast$). Thus, by using the Goethals--Seidel 
array, they give four skew Hadamard matrices of order 188.
We have now constructed an SDS with parameters $(47;23,21,19,19;35)$
and type (ks$\ast$$\ast$). It gives a new skew Hadamard matrix
of order 188. Here is this SDS:
\begin{eqnarray*}
&& \{1,2,3,4,6,7,8,9,12,14,16,17,18,21,24,25,27,28,32,34,36,37,42\}, \\
&& \{0,6,8,10,11,14,17,18,19,21,23,24,26,28,29,30,33,36,37,39,41\}, \\
&& \{0,1,2,5,6,8,9,15,16,19,21,23,27,28,33,36,38,39,40\}, \\
&& \{0,2,3,4,7,8,9,10,12,18,21,23,24,25,26,30,34,35,44\}.
\end{eqnarray*}
The first set is the $(47;23;11)$ skew difference set consisting 
of all nonzero squares in $\bZ_{47}$.

\subsection{Multicirculant Williamson matrices of order $49$} 
According to \cite{KS} there are no cyclic SDSs of type (ssss) with the
parameters $(7^2;21,21,21,21;35)$. On the other hand, elementary 
abelian SDSs having the same type and parameters exist;
an example due to R.M. Wilson is given in \cite{QX}.
We have constructed another such SDS, not equivalent to Wilson's example.

Let $F_{49}=\bZ_7[x]/(x^2-3)$ and let us identify $x$ with its image 
in $F_{49}$. Our SDS consists of four symmetric blocks
$A_i=A'_i\cup \{0\} \cup(-A'_i)$ in $\pA=(F_{49},+)$, where
\begin{eqnarray*}
&& A'_1 = \{1,2,2x,x+2,2x+2,2x-1,3x+1,3x-2,3x\pm3\}, \\
&& A'_2 = \{2,x,2x,x+2,x-1,x\pm3,2x+3,3x+2,3x+3\}, \\
&& A'_3 = \{2,2x,x+1,x+3,x-2,2x-2,2x\pm3,3x-2,3x-3\}, \\
&& A'_4 = \{3,3x,x+2,x+3,2x+1,2x+3,3x-1,3x-2,3x\pm3\}.
\end{eqnarray*}

Note that $\{\pm1,\pm3x\}$ is the unique subgroup of order 4 
in $F_{49}^*$ and none of the $A_i$ contains this subgroup.
In Wilson's example one of the four blocks contains the subgroup
of order 4. By using this fact and an argument from \ref{Slucaj27},
it is easy to show that the two examples are not equivalent.
Both examples give rise to multicirculant 
Williamson matrices of order 49.

\subsection{A new skew Hadamard matrix of order $4\cdot61$}
We have constructed a cyclic SDS $(A_i)$ with parameters
$(61;30,28,27,24;48)$ and symmetry type (k$\ast$$\ast$s).
The four blocks are:
\begin{eqnarray*}
A_1 &=& \{ 1,6,7,9,13,16,17,18,20,22,24,25,27,28,30,32,35,38, \\
&& 40,42,46,47,49,50,51,53,56,57,58,59 \}, \\
A_2 &=& \{ 0,1,2,3,7,11,12,13,14,15,19,21,22,24,26,28,29,30, \\
&& 33,34,35,39,42,47,48,58,59,60 \}, \\
A_3 &=& \{ 2,3,4,5,11,16,19,20,21,22,25,26,27,29,32,33,36,39, \\
&& 40,41,42,45,46,49,50,52,58 \}, \\
A_4 &=& \{ 7,8,10,12,15,16,18,20,24,25,27,30,31,34,36,37,41, \\
&& 43,45,46,49,51,53,54 \}. 
\end{eqnarray*}
By using the Goethals--Seidel array, we obtain a new skew Hadamard
matrix of order $4\cdot61$.

\subsection{A new skew Hadamard matrix of order $4\cdot127$}
We have constructed a cyclic SDS $(A_i)$ with parameters
$(127;57,57,57;76)$ and symmetry type (ks$\ast$$\ast$).
As 127 is a prime, we have $\pA=(\bZ_{127},+)$.
Let $H=\{ 1,2,4,8,16,32,64 \}$, the subgroup of $\bZ_{127}^*$
of order 7. We enumerate its 18 cosets as $\al_i$, $0\le i\le 17$,
such that $\al_{2i+1}={-1}\cdot\al_{2i}$, $0\le i\le 8$.
For even indexes we have
\[
\begin{array}{lllll}
\al_0=H, \quad & \al_2=3H, \quad & \al_4=5H, \quad & \al_6=7H, \quad & \al_8=9H, \\
\al_{10}=11H, & \al_{12}=13H, & \al_{14}=19H, & \al_{16}=21H. & 
\end{array}
\]
We use the index sets:
\begin{eqnarray*}
J_1 &=& \{ 0,1,2,3,6,7,16,17 \}, \\
J_2 &=& \{ 4,6,7,11,13,14,15,16 \}, \\
J_3 &=& \{ 0,4,5,7,11,12,15,16 \}
\end{eqnarray*}
to define the three blocks by
\[ A_i = \{ 0 \} \cup \bigcup_{k\in J_i} \al_k, \quad 1\le i\le 3. \]
By combining this SDS with the classical Paley skew $(127;63;31)$ difference 
set, we obtain an SDS with parameters $(127;63,57,57,57;107)$ and
type (ks$\ast$$\ast$). By using the Goethals--Seidel array, it gives 
a new skew Hadamard matrix of order $4\cdot127$.

Note that the above SDS $(A_1,A_2,A_3)$ is a difference family and so
it gives a balanced incomplete block design (BIBD) with parameters
$(v,k,\la)=(127,57,76)$.


\begin{thebibliography}{99}
\bibitem{RV}
Rudy van Vliet, The exhaustive search for Williamson matrices of
order 45 and 51 is attributed to him, see e.g. \cite{HKT}.
\end{thebibliography}
\end{document}